\documentclass[1pt,notitlepage,twoside,a4paper]{amsart}

\usepackage{amsmath,amssymb,enumerate}

\usepackage{epsfig,fancyhdr,color}%,showkeys,amsmidx

\usepackage{amssymb}
\usepackage{amsmath,amsthm}    
\usepackage{latexsym}
\usepackage{amscd} 
\usepackage{psfrag}
\usepackage{graphicx} 
\usepackage[latin1]{inputenc}  
\usepackage[all]{xy} 
\usepackage[mathcal]{eucal}

\theoremstyle{definition}

\theoremstyle{remark}

\def\interieur#1{\mathord{\mathop{\kern 0pt #1}\limits^\circ}}

\definecolor{NoteColor}{rgb}{1,0,0}

\title[Nicolas-Auguste Tissot]{Nicolas-Auguste Tissot: 
A link  between cartography and quasiconformal theory}

\author{Athanase Papadopoulos}
\address{A. Papadopoulos, Institut de Recherche Math{\'e}matique Avanc\'ee,
Universit{\'e} de Strasbourg and CNRS,
7 rue Ren\'e Descartes,
 67084 Strasbourg Cedex, France}

 \date{\today}

\begin{document}
\begin{abstract}  
Nicolas-Auguste Tissot (1824--1897) published a series of papers on cartography in which he introduced a tool which became known later on, among geographers, under the name of the \emph{Tissot indicatrix}. This tool was broadly used during the twentieth century in the theory and in the practical aspects of the drawing of geographical maps. The Tissot indicatrix is a graphical representation of a field of ellipses on a map that describes its distortion. Tissot studied extensively, from a mathematical viewpoint, the distortion of mappings from the sphere onto the Euclidean plane that are used in drawing geographical maps, and more generally he developed a theory for the distorsion of mappings between general surfaces. His ideas are at the heart of the work on quasiconformal mappings that was developed several decades after him by Gr\"otzsch, Lavrentieff, Ahlfors and Teichm\"uller. 
Gr\"otzsch mentions the work of Tissot and he uses the terminology related to his name (in particular, Gr\"otzsch uses the Tissot indicatrix).
 Teichm\"uller mentions the name of Tissot in a historical section in one of his fundamental papers where he claims that quasiconformal mappings were used by geographers, but without giving any hint about the nature of Tissot's work. The name of Tissot is also missing from all the historical surveys on quasiconformal mappings. In the present paper, we report on this work of Tissot. We shall also mention some related works on cartography, on the differential geometry of surfaces,  and on the theory of quasiconformal mappings. This will place  Tissot's work in its proper context.

\bigskip



\end{abstract}
  \maketitle

\tableofcontents

\section{Introduction: From geography to quasiconformal mappings}
 
The theory of quasiconformal mappings (like that of conformal mappings) can be traced back to old geography.  Geographers, since Greek antiquity, were interested in the question of drawing maps that minimize distortion. Here, the word ``distortion" may refer to angle, area, or length, or a combination of these notions. The early geographers knew that a mapping from the sphere (or a spheroid)\footnote{The claim that the Earth is spheroidal and not spherical, more precisely, that it is slightly flattened at the poles, was first made by Newton, who concluded in his \emph{Principia} that this flatness, which he expected to be of the order of 1/230,  is due to the Earth's rotation. This was confirmed by several expeditions in the eighteenth century, which included especially French scientists, whose aim was to make precise measurements of the meridians near the poles.} to the Euclidean plane cannot preserve these three parameters at the same time, and they searched for a ``best compromise."  Making  precise such a notion of compromise between the three distortions we mentioned depended on the practical use for which the map was intended. 

Several mathematicians -- some of them among the most prominent -- were also geographers. One can mention Ptolemy, Euler, Lambert, Lagrange, and Gauss, and there are many others. Most of them were interested at the same time in the mathematical theory and the practical art of map drawing. Euler, at the Academy of Sciences of Saint Petersburg, besides being a mathematician, had the official status of a cartographer, and was one of the leaders of the huge project of drawing maps of the new Russian Empire. Motivated by cartography, he was led to study questions regarding conformal and almost-conformal mappings. In 1777, he published three memoirs on mappings from the sphere to the Euclidean plane: \emph{De repraesentatione superficiei sphaericae super plano} (On the representation of  spherical surfaces on a plane) \cite{Euler-rep-1777},  \emph{De proiectione geographica superficiei sphaericae} (On the geographical projections of spherical surfaces) \cite{Euler-pro-1777}  and  \emph{De proiectione geographica Deslisliana in mappa generali imperii russici usitata} (On Delisle's geographic projection used in the general map of the Russian empire)  \cite{Euler-pro-Desli-1777}. In the first memoir, \cite{Euler-rep-1777}, Euler examines several projections of the sphere, searching systematically for the partial differential equations that they satisfy. He recalls that there is no ``perfect" mapping from the sphere onto a plane, and he considers the question of finding mappings with the least distortion. He  highlights several kinds of properties that one may naturally ask for geographical maps: conformality; sending meridians to curves normal to a given axis and parallels to lines parallel to that axis; preservation of area (up to scale), etc. Most of all, Euler considered, like Lambert at about the same time, arbitrary mappings from the sphere onto a flat surface, and not only mappings obtained as a projection from a given point onto a plane, a cylinder or a cone, etc.

Lagrange, in his study of geographical maps (see \cite{Lagrange1779} and \cite{Lagrange1779a}), also considered general mappings from the sphere onto a flat surface. Like Euler and Lambert before him, he characterized angle-preserving general mappings analytically, and he applied his methods to the particular case of mappings from the sphere onto the plane that send meridians and parallels to circles.  He writes that ``these are the only curves which one can easily employ in the construction of geographic maps" (\cite{Lagrange1779} pp. 642). Lagrange introduced a function called the \emph{magnification ratio}. The name given to this function clearly indicates its meaning. For conformal mappings, the magnification ratio depends only on the point on the sphere (and not on the direction at that point). Lagrange gave formulae for this parameter for general conformal mappings. We note right away that this notion played an important role in the later work of Chebyshev on the drawing of geographical maps: the latter noticed that the magnification ratio, under the condition that the distortion of the map is minimal in some appropriate sense, satisfies Laplace's equation, and he showed that under some natural conditions, the best angle-preserving geographical maps of a certain region having a reasonable size are those for which the magnification ratio function is constant on the frontier of the region. Chebyshev thus reduced the question of finding the best  geographical map to a problem in potential theory, namely, solving the Laplace equation on a certain domain of the sphere for a function which is constant on the boundary of the domain. We refer the reader to the papers \cite{Cheb1} and \cite{Cheb2} by Chebyshev (also included in his \emph{Collected works} \cite{T-oeuvres}), and to the recent survey \cite{Papa-Chebyshev} for some details on the subject.

It is also worth recalling that the motivation of several prominent mathematicians to study the differential geometry of surfaces came from geography. A famous case is Gauss, who declares explicitly in the preface of his paper \cite{Gauss-Copenhagen} that his aim is only to construct geographical maps and to study the general principles of geodesy for the task of land surveying. Gauss's paper \cite{Gauss-Copenhagen} is titled \emph{Allgemeine Aufl\"osung der Aufgabe, die Teile einer gegebenen Fl\"ache auf einer andern gegebenen Fl\"ache so abzubilden dass die Abbildung dem Abgebildeten in den kleinisten Theilen \"ahnlich wird.} (General solution of the problem: to represent the parts of a given surface on another so that  the smallest parts of the representation shall be similar to the corresponding parts of the surface represented). This is an epoch-making paper, in which the author proves that every sufficiently small neighborhood  of any point in an arbitrary real-analytic surface can be mapped conformally onto a subset of the Euclidean plane.  

There are many other examples of mathematical works on the geometry of surfaces motivated by geography. We refer the reader to the recent survey \cite{Papa-qc}.

After differential geometry, we arrive at quasiconformal mappings.  Teichm\"uller's paper \cite{T20} contains a section on the origin of these mappings in which he mentions the name of the French geographer Tissot. Gr\"otzsch, in his paper \cite{Groetzsch1930} mentions several times the name of Tissot and refers to his work. But the name of Tissot is never mentioned in any paper on the history of quasiconformal mappings, and it seems that the references to him by Gr\"otzsch and by Teichm\"uller's remained unnoticed. Yet the work of Tissot, which is well known among geographers, and in particular in the drawing of geographical maps, is very closely related to quasiconformal mappings. Let us mention right away some of Tissot's important results.
   \begin{enumerate}
\item  Tissot invented a device which provides a visual measure for the distortion of a geographical map. This is a field of infinitesimal ellipses which are the image of a field of infinitesimal circles on the region on the sphere which is represented. This field of infinitesimal ellipses is characterized by two parameters: their relative size, and the ratio of their two axes (the major axis divided by the minor axis).\footnote{The expression ``infinitesimal circle" means here, as is usual in the theory of quasiconformal mappings, a circle on the tangent space at a point. In practice it is a circle on the surface which has a ``tiny radius." These circles, on the domain surfaces, are all supposed to have the same small size, so that the collection of relative sizes of the image ellipses is a meaningful object.}
 These parameters represent respectively the local area distortion and the angle distortion of the mapping that is used to draw the geographical map.  
 This device is known  among geographers as the \emph{Tissot indicatrix}. An example is given in Figure \ref{Tissot-Cassini}.
  From the differential geometric point of view, the Tissot indicatrix gives an information on the metric tensor of the metric obtained by pushing forward the metric of the sphere or the spheroid by the projection used.
  \begin{figure}[htbp]
\centering
\includegraphics[width=10cm]{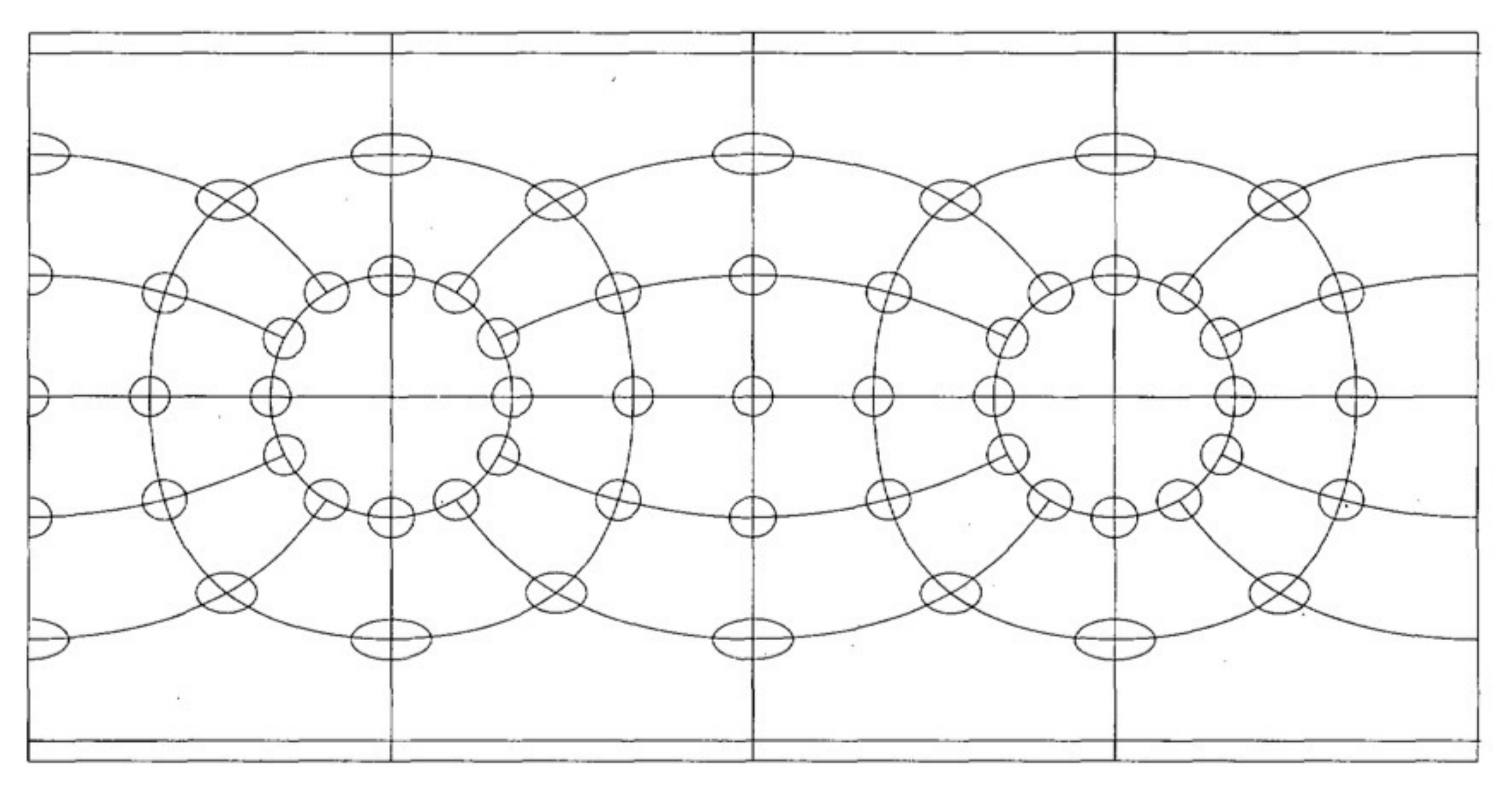}
\caption{\small{The Tissot indicatrix for the Cassini projection; from the \emph{Album of map projections} \cite{Snyder-Album} pp. 26.}} \label{Tissot-Cassini}
\end{figure}

\item Tissot noticed that for a given mapping between two surfaces, there is, at each point of the domain, a pair of orthogonal directions which are sent to a pair of orthogonal directions on the image surface. Unless the mapping is angle-preserving at some point, the pair of orthogonal directions is unique. The orthogonal directions at the various points on the two surfaces define a pair of orthogonal foliations which are preserved by the mapping. We shall comment later on the uniqueness of the two fields of directions that are referred to in the above result.

An important observation made by Tissot right at the beginning of his memoir \cite{Tissot1881} (p. 1) is that  
the most appropriate mode of projection in order to build a geographical map depends on the shape of the region, that is, on the shape of its boundary as a subset of the sphere (or the spheroid). Tissot discovered that in order for the map to minimize the appropriately defined distortion, a certain function $\lambda$, defined by setting 
\[d\sigma^2=(1+\lambda)^2ds^2,\]
must be minimized in some sense.
Here, $ds$ and $d\sigma$ are the line elements at the source and the target surfaces respectively. The minimality of $\lambda$ may mean, for example, that the value of the gradient of its square must be the smallest possible. Darboux later on connected this work of Tissot with the work of Chebyshev that we mentioned above. We shall comment on this fact later in this paper.
 
 \item Tissot constructed mappings that realize the required properties of minimal distortion. 
 
 \end{enumerate}
 
 Most of the mathematicians who preceded Tissot (Lagrange,  Euler, etc.), and Chebyshev who was his contemporary, studied conformal maps between surfaces, with some distortion parameter to be minimized. In contrast, Tissot studied extensively distortion-minimizing \emph{non-conformal} maps. 

  We shall elaborate on Tissot's ideas in \S\,\ref{s:cart} below.

 Beyond the work of Tissot, a theory of ``almost-conformal" mappings whose main object was to find mappings with least distortion between two surfaces had been inherent in the works of geographers for more than two thousand years. 
Geographical maps were used for various practical needs, including navigation,  surveying for the purpose of collecting taxes, etc., 
 but the theory behind these maps was an important element in the exact sciences.

   It is always interesting how very simple ideas lead to elaborate theories.
   
\section{On the works of Nicolas-Auguste Tissot}

Nicolas-Auguste Tissot (1824--1897)  was born (like Poincar\'e whom we shall mention very soon) in Nancy. In an article by A.-P.  d'Avezac published in the \emph{Bulletin de la Soci\'et\'e de G\'eographie} \cite{Avezac}, we read the following (p. 127)\footnote{This quote and the next one are my translations from the French.}:

\begin{quote}\small 
Mr. Nicolas-Auguste Tissot, an ancient \emph{capitaine du g\'enie}, and today assistant professor (\emph{r\'ep\'etiteur}) in geodesy at the  \'Ecole Polytechnique, presented to the Academy of Sciences a series of communications on geographical maps, whose outline is reproduced in the \emph{Comptes Rendus 
des S\'eances} (Novembre 7, 1859, 49, pp. 673--676; March 5, 1860, pp. 474--476; and December 17, 1860, pp. 964--969). Their main goal is to determine the rule according to which the deformation is produced around each point, regardless of the system of representation;  and to compare, from this point of view, the various systems that are used and proposed in the drawing of mappae mundi, and finally, to find, in this way, the best mode of projection for each particular country. In the investigation of this particular question, he outstripped Mr. Airy.\footnote{The reference is to the British astronomer George Biddell Airy (1801--1892).}
\end{quote}

In \S\,\ref{s:cart}, we shall describe the content of the \emph{Comptes Rendus} notes referred to by d'Avezac. 

 Tissot defended his doctoral thesis on November 17, 1851, cf. \cite{Tissot-these}. On the cover page of the thesis, the author is described as ``M. A. Tissot, Ex-Capitaine du G\'enie." The work presented contains in fact two theses,\footnote{In France, the \emph{Doctorat d'\'Etat} consisted of two theses; the second one was not necessarily original,  and usually consisted in  the survey of some work in a field which is different from the candidate's field of expertise. The subject of this thesis was chosen by the faculty and communicated to the PhD student a few weeks before the defence. The Doctorat d'\'Etat disappeared in the 1990s.} the first one on Mechanics, titled \emph{Mouvement d'un point mat\'eriel pesant sur la sph\`ere} (Motion of a material point moving on the sphere) and the second one on astronomy, titled \emph{Sur la d\'etermination des orbites des plan\`etes et des com\`etes} (On the determination of the orbits of planets and comets). In his first thesis, Tissot studies the motion of a point on a sphere under the action of gravity. In the second thesis, Tissot describes various methods of determining distances between celestial bodies making use of works of Cauchy and of Lagrange on the determination of the orbits of celestial bodies.
   The techniques used in the two memoirs are those of differential geometry, differential equations and approximation theory. They are the same techniques that Tissot used later on in his work on cartography.
 
On the title page of Tissot's memoir \cite{Tissot1881} (1881), the expression \emph{Examinateur \`a l'\'Ecole Polytechnique} follows his name. In his \emph{\'Eloge historique de Henri Poincar\'e}  \cite{Darboux-Eloge}, Darboux  relates the following episode about Tissot, examining Poincar\'e at the entrance exam:\footnote{Poincar\'e entered the \'Ecole Polytechnique in 1873. It is not unusual that at such an examination, some students listen to the examination of other students. In fact, a student is given a question or a set of questions which he is asked to prepare while another student (who had already been given some time to prepare his questions) is explaining his solutions at the blackboard, in the same room.}
 \begin{quote}\small
 Before questioning Poincar\'e, Mr. Tissot suspended the exam during 45 minutes: we thought it was the time he needed to prepare a sophisticated question. Mr. Tissot came back with a question of the Second Book of Geometry. Poincar\'e drew a formless circle, he marked the lines and the points indicated by the examiner, then, after wandering long enough in front of the blackboard, with his eyes fixed on the ground, he concluded loudly: ``It all comes down to proving the equality $AB=CD$. This is a consequence of the theory of mutual polars, applied to the two lines." Mr. Tissot interrupted him:  ``Very good, Sir, but I want a more elementary solution." Poincar\'e started wandering again, this time not in front of the blackboard, but in front of the table of the examiner, facing him, almost unconscious of his acts; then suddenly he developed a trigonometric solution. 
Mr. Tissot objected: ``I would like you to stay in elementary geometry." Almost immediately after that, the examiner of elementary geometry was given satisfaction. He warmly congratulated the examinee and announced that he deserves the highest grade.
  \end{quote}

Besides his work on geographical maps, Tissot wrote several papers on elementary geometry; see e.g. \cite{Tissot1}, \cite{Tissot2} and \cite{Tissot3}. We mention in this respect that several prominent French mathematicians of the nineteenth and the beginning of the twentieth century were very much interested in such problems. The names include  Serret, Catalan, Laguerre, Darboux, Hadamard and Lebesgue, and there are others.

\section{Tissot's work on cartography} \label{s:cart}

In the years 1856 and 1858, Tissot published several papers and \emph{Comptes Rendus} notes on cartography  in which he made a detailed analysis of the distortion of some known geographical maps (see \cite{Tissot-CR00}, \cite{Tissot-1858}, \cite{Tissot-CR0}). He started developing his own theory in three \emph{Comptes Rendus} notes \cite{Tissot-CR1}, \cite{Tissot-CR2} and \cite{Tissot-CR3}, written in the years 1859--1860.

In \cite{Tissot-CR1}, Tissot starts with the following observation. To construct a geographical map, or to represent one surface on another, one first chooses on each of the two surfaces a decomposition into infinitesimal parallelograms using two systems of curves (in reality, these two ``systems of curves" are two transverse foliations), in such a way that the lines of the first system are sent to the lines on the second.  In this way, points and their images are encoded as intersection points of corresponding lines.

Encoding points as intersections of two lines was widely used in projective geometry. At the practical level, it was used by engineering students in France in their study of  ``dessin industriel", by students in arts who were asked to magnify a drawing using graph paper, etc.

Based on the existence of this decomposition into nets of infinitesimal parallelograms, Tissot states the following:

\medskip

\emph{Any representation from a surface onto another may be replaced, at each point, by an orthogonal projection made at an appropriate scale.}

\medskip

The proof is given in his later papers. Tissot then gives in \cite{Tissot-CR1} the following principles that concern angle, distance, and area deformation.
These principles play a very important role in the rest of his work.

\medskip

\emph{1.--- For any kind of representation of a surface onto another, there exists, at every point of the first surface, two perpendicular tangents, which are unique unless the angles are preserved at that point, such that the images of these two tangents are perpendicular on the second surface.}

\medskip

Tissot calls the tangents referred to in this statement \emph{principal tangents}.

\medskip

\emph{2.--- The directions of the principal tangents are those at which the ratio of lengths of the corresponding infinitesimal elements attains its greatest and smallest values.}

\medskip 
Tissot denotes these two values by $a$ and $b$, $a>b$.

\medskip

\emph{3.--- In order to find the image of an infinitely small figure drawn in the tangent plane of the first surface, one uses the following procedure: Superpose the two tangent planes in such a way that the principal tangents coincide. Turn the first tangent plane by an angle whose cosine is $b/a$. Use the orthogonal projection given by the two decompositions by parallelograms on the two surfaces, then modify the image of the figure by applying in the second tangent plane a magnification of ratio $a$.}

\medskip

Tissot notes that the image of an infinitely small figure around a point on the first surface which is a circle of radius 1 centered at the given point is an ellipse whose major and minor axes are $a$ and $b$.
He then outlines a practical way to find the major and minor axes of these ellipses, and he provides formulae for them.  This is the basis of the theory of what became known later as the \emph{Tissot indicatrix}.

In the note \cite{Tissot-CR2}, Tissot says that at any point of the Terrestrial globe,  the angle  which is deformed most by a geographical map is never the angle made by the meridian and the parallel at that point. Likewise, the directions which are stretched most, or contracted most, do not coincide with those of the meridian and the parallel unless the images of these directions are at right angles. He says that using the rules he gave in the previous note,  one can easily compute the greatest alterations of angles and distances. This is the basis of the theory which is behind the maps that he will describe later on, in which he seeks for the most advantageous properties for what concerns distortion. 

Tissot declares that the projection that he shall adopt depends on the particular country that has to be represented, in particular, its position with respect to the equator, its size, and the form of its contour. He then makes a list of 12 known projections or families of projections that are known, for which he gives formulae for the angle, distance, and area deformations. He compares the usefulness of these projections. The central question is again that of finding the best  possible projection.

The third note, \cite{Tissot-CR3}, is more technical than the first two. Tissot mentions there the difficulties that arise in applying the theory to a country with a great area, like Russia. He gives several formulae for the projections of countries with a reasonable surface area, like France and Spain.
These formulae will be discussed more thoroughly in his long memoir \cite{Tissot1881} published several years later. In a note on pp. 2 of this memoir, Tissot declares that after he published the first of his \emph{Comptes Rendus} notes on the subject, the statements that he gave there without proof were reproduced by A. Germain in his \emph{Trait\'e des projections des cartes g\'eographiques} \cite{Germain} and by U. Dini in his memoir \emph{Sopra alcuni punti della teoria delle superfici} \cite{Dini}. He says that Germain and Dini gave their own proofs of these statements, which are nevertheless more complicated than those he had in mind and which he gives in the present memoir.
He also declares that Dini showed that the whole theory of curvature of surfaces may be deduced from the general theory that Tissot developed. In fact, Dini applied this theory to the representation of a surface on a sphere, using Gauss's methods. Tissot also declares that his ideas were used in astronomy, by Herv\'e Faye, in his \emph{Cours d'astronomie de l'\'Ecole polytechnique} \cite{Faye}.  The texts of the two \emph{Comptes rendus} notes \cite{Tissot-CR3} and \cite{Tissot-CR0} of Tissot are reproduced in the treatise of Germain \cite{Germain}.

 Tissot developed his complete theory several years after his three \emph{Comptes Rendus} notes, first in installments, in 1878--1880,  (cf.  Tissot \cite{Tissot1878}--\cite{Tissot1880}), and then in the long memoir \cite{Tissot1881}, which contains the work in \cite{Tissot1878}--\cite{Tissot1880} together with additional material.

 The memoir \cite{Tissot1881} is 337 pages long. It consists of a preamble and 4 chapters. We briefly review its content.

In the first chapter, Tissot studies the general properties of surfaces and the distortions of maps between them. He provides the proof of Statement 1 which we quoted above, asserting d the existence and uniqueness of the pairs of field of perpendicular lines that are sent to each other by the maps. He also gives the following complement:

\medskip
\emph{We can vary in all possible manners the fields of perpendicular lines,  obtaining  an infinite number of decompositions into infinitesimal rectangles, each leading to a map between surfaces and having the property that these fields are the unique invariant orthogonal fields. 
}

\medskip

The first chapter contains, besides the proofs of the results announced in the \emph{Comptes Rendus} notes, a description of the Tissot indicatrix, a detailed study of its properties, and a method of calculating the various deformations of angle, length and area using this indicatrix. 

Practically, on a given geographical map, the Tissot indicatrix is a field of ellipses drawn at certain points (often located at intersections of the images of the parallels and the meridians that are drawn on the map). These ellipses illustrate the distortion at these points. See e.g. the maps drawn in the \emph{Album of map projections} \cite{Snyder-Album}.

In terms of the differential geometry of surfaces, the Tissot indicatrix is a representation of the quadratic form that defines the metric tensor at each point of the geographic map. In the case of a conformal map, the Tissot indicatrix is simply a field of circles. 
 
 Tissot explains how one can draw the axes of the indicatrix. The two fileds of axes, when they are integrated over the surface, give the two orthogonal foliations, forming what which Tissot calls \emph{canevas}. The development of the theory is an ingenious combination of elementary geometry, infinitesimal calculus, and differential equations. We recall that the results apply to maps between  differentiable surfaces embedded in space, and not only between the sphere and the Euclidean plane. Tissot studies in particular the case of surfaces of revolution. 

 In the section starting at pp. 21 of the memoir, Tissot gives the formulae for the major and minor axes of the ellipses and its inclination for general maps between surfaces, in terms of the local parametrisations of these surfaces.

Tissot starts with two surfaces defined using orthogonal coordinate systems $(x,y,z)$ and $(x',y',z')$ respectively, functions of two variables $l$ and $m$. The coordinates are chosen so that the mapping between the surfaces is defined by equating the six coordinates with appropriate functions of $l$ and $m$. He sets
\[L=\left[\left(\frac{dx}{dl}\right)^2+\left(\frac{dy}{dl}\right)^2+\left(\frac{dz}{dl}\right)^2\right]^{\frac{1}{2}}
,\]
\[L'=\left[\left(\frac{dx'}{dl}\right)^2+\left(\frac{dy'}{dl}\right)^2+\left(\frac{dz'}{dl}\right)^2\right]^{\frac{1}{2}}
,\]
\[M=\left[\left(\frac{dx}{dm}\right)^2+\left(\frac{dy}{dm}\right)^2+\left(\frac{dz}{dm}\right)^2\right]^{\frac{1}{2}}
,\]
\[M'=\left[\left(\frac{dx'}{dm}\right)^2+\left(\frac{dy'}{dm}\right)^2+\left(\frac{dz'}{dm}\right)^2\right]^{\frac{1}{2}}
.\]
The lengths of two infinitesimal parallelograms on the two canevas that correspond to each other will be, for the first surface, $Ldl$ and $Mdm$, and on the second one, $L'dl$ and $M'dm$. 
Setting $h$ and $k$ to be the ratios of the sides of these parallelograms and $\Theta$ and $\Theta'$ their angles, we have
\[h=\frac{L}{L'},\]
\[k=\frac{M}{M'},\]
\[\cos \Theta=\frac{1}{LM}\left[\frac{dx}{dl}\frac{dx}{dm}+\frac{dy}{dl}\frac{dy}{dm}+\frac{dz}{dl}\frac{dz}{dm}
\right]
\] 
and
\[\cos \Theta=\frac{1}{L'M'}\left[\frac{dx'}{dl}\frac{dx'}{dm}+\frac{dy'}{dl}\frac{dy'}{dm}+\frac{dz'}{dl}\frac{dz'}{dm}
\right].
\] 
Depending on the practical applications needed, the formulae for the major and minor axes and the inclination of the ellipses defining the Tissot indicatrix are derived from these formulae. We refer the reader to pages 22ff. of Tissot's  memoir, where he considers a variety of special cases in which he derives precise formulae. For instance, in the case where the mapping is such that the ratios of lengths in the direction of the two linear elements have the same common value $h$, then he finds the following formulae for the major and minor axes of the ellipses:
\[b=h\frac{\displaystyle\frac{\cos\Theta'}{2}}{\displaystyle\frac{\cos\Theta}{2}}\]
and 
\[b=h\frac{\displaystyle\frac{\sin\Theta'}{2}}{\displaystyle\frac{\sin\Theta}{2}}\]
Many other cases are considered.

Tissot uses extensively power series expansion. We shall see below that Darboux, in Part II of his paper \cite{Darboux-Tissot} on the work of Tissot, develops the same theory in the setting of Gauss's intrinsic geometry of surfaces. It is interesting to note that after he explains the methods of Tissot, Darboux says: ``We can approach the same theory in a more rigorous way, where the power series developments intervene only at the end, and in more general conditions."

 The second chapter concerns the practical applications of the theory developed in Chapter 1. It is dedicated to the answer to the question: 
 
 \emph{How to find the projection which is most appropriate to a given country?} 
 
 In the introduction of this chapter, Tissot says that in the drawing of maps that are intended for the use in public services, and in particular in the services of the army, the most important question is the reproduction of angles. He says that for the map to be useful, as a topographical tool, the distortion of angles must be very small. He recalls that the distortion of distances varies from point to point, and he says that the \emph{supremum} of this distortion must be reduced to a minimum. Another rule he mentions is that in working with geographical maps,  \emph{formulae must be simple}. He suggests equipping the map with a network of curves, using colors, which are loci where the distortion of length is constant. He shows that these curves will be in general algebraic of degree two, and in most cases they will be ellipses. Furthermore, for a given country that has to be represented, a point has to be chosen, called the \emph{central point}, and the position of this point must be determined according to the required maximizing properties. The meridian and the parallel passing through this point will be called the \emph{mean meridian} and the \emph{mean parallel}.  This chapter contains many formulae,  concerning types of various projections. 
As an illustration, let us mention one of the formulae (p. 71), which is also contained in his \emph{Comptes Rendus} note \cite{Tissot-CR3} written 21 years before. It involves some notation.

Let $l$ and $m$ denote the latitude and longitude of a point on the surface of the Earth, $r$ the radius of the terrestrial parallel at the latitude $l$ and $s$ the arc of meridian contained between the mean parallel and the parallel at latitude $l$. If $x$ and $y$ denote the coordinates of the corresponding point on the map with respect to a perpendicular system of axes, then the projection is given by the following formulae
\[ x=s+\frac{1}{2}rm^2\sin l
\]
and 
\[y=rm(1+\frac{1}{6}m^2\cos 2l)
.\]
The other formulae contain higher order terms.

Then follows a long development in which Tissot  simplifies the formulae using power series expansions, with practical applications concerning the drawing of maps of France, Spain and other countries.

The last two chapters of the memoir \cite{Tissot1881} are more computational than the first two. Tissot compares there the distortions of various projections.

Let us also mention the so-called \emph{Tissot projection} used in cartography. It is described in Chapter 2 of the memoir \cite{Tissot1881}. From the point of view of minimizing distortion, this map makes a compromise between conformal  and equal area projection. The mathematical tools behind  the definition are based again on power series expansions.   The Tissot projection was used by the French army cartography department; see \cite{Snyder-F}. It was also used in the cadastral survey of Egypt; cf. \cite{Lyons}. In fact, Tissot mentions in his memoir \cite{Tissot1881} pp. 71, that his projection can be used for  a country situated  in the lune between two meridians which are not very distant from one another, and as an example he gives precisely Egypt.

Darboux was interested in the work of Tissot on geography, and in particular, in his projection described in Chapter 2 of his memoir \cite{Tissot1881}. He wrote a paper on Tissot's work, \cite{Darboux-Tissot}.  In the introduction of that paper, Darboux recalls that using power series expansions, Tissot gave a new method for the representation of a given country. The region that has to be represented has to be small enough compared to the Earth. He also recalls that the series expansions reduce to a minimum amount the distortions of angles and distances, but he says that ``[Tissot's] exposition appeared to me a little bit confused, and it seems to me that while we can stay in the same vein, we can follow the following method." He then provides another method to construct Tissot's projection.  
 In this method, using a system of principal tangent lines, the linear element at a point of the surface is written as
   \begin{equation}\label{D:1}
   ds^2=dx^2+dy^2+\left(\frac{xdx}{R}+\frac{ydy}{R'}\right)^2,   
   \end{equation}
   where $R$ and $R'$ are the principal radii  of curvature at the given point. In this formula, the terms of the third order and more in  $x$ and $y$ have been neglected.
   
   Then, one looks for power series,
   \[\alpha=f(x,y)\]
   and
   \[\beta=\phi(x,y)\]
   such that 
   \[d\alpha^2+d\beta^2\]
   coincides with $ds^2$ for terms up to a certain order (the highest possible). 
   
   The condition that the map is a similarity at that point leads to an equation of the form
      \[d\alpha^2+d\beta^2=(1+\lambda)^2ds^2\]
      where $\lambda$ is a function.
      
     Setting $\lambda_0$ to be the second order homogeneous component of  $\lambda$, Darboux is led to the equation
      \[\lambda_0=\frac{x^2+y^2}{4RR'}+A(x^2-y^2)+2Bxy\]
      where $A$ and $B$ are constants.
      
 Darboux   then says that Tissot determines $A$ and $B$ by an ingenious device. 
 Let us admit for simplicity that the Earth is spherical, that is, let us take $R=R'=1$. 
 Then Tissot's method amounts to searching, among all the conics with equation
 \[\lambda =\frac{\alpha^2+\beta^2}{4}+A(\alpha^2-\beta^2)+2B\alpha\beta=\mathrm{const.},\]
 the one which fits better the form of the surface to be represented, and at the same time corresponds to the smallest value of $\lambda$, with some considerations to be taken into account, for instance, the mean value of the square of the gradient of $\lambda$ is sought to be smallest possible.
 
 Darboux says that in the case where the curve is an ellipse, one obtains a ``remarkable result,", namely, that \emph{the value of $\lambda$ is maximal on the boundary of the ellipse, and the value of $\lambda$ on this boundary is constant}. Darboux (like Tissot before him)  gives the value of this constant. 
Darboux notes that there is a relation with Chebyshev's  theorem that we mentioned above, and he makes this relation explicit.  In both theories, an integral has to be minimized, and if one takes the same degree of approximation, the two integrals are the same.
Darboux then explains that conversely, the theorem of Chebyshev justifies the ``rather elementary" reasoning of Tissot and his use of the set of conics. Indeed, by the theorem of Chebyshev, $\lambda$ has to be constant on the boundary of the region to be represented and since, with the approximation made, it has to be constant on concentric conics, it is natural to choose among these conics the one which covers best the given region or, more precisely, its projection on the tangent plane, while having the smallest $\lambda$.

   After that, Darboux extends the whole theory of Tissot to maps between surfaces in the setting of Gauss's theory, starting with  arbitrary curvilinear coordinates with length element
   \[ds^2=Edu^2+2Fdudv+Gdv^2\]
   and using the theory of conformal representations.
   He makes the relation between the work of Tissot and the works of Gauss, Tchebyshev and Beltrami.

Since we mentioned the works  of Lagrange and Chebyshev on cartography, let us note that Darboux wrote two papers which are directly motivated by these works.  In the first paper, \cite{Darboux-Lagrange}, Darboux gives a proof of a problem addressed by Lagrange in   the paper \cite{Lagrange1779} on cartography that we already mentioned. The question concerns a constant which Lagrange calls the ``exponent of the projection." This question is reduced to the following problem in elementary geometry: \emph{Given three points on the sphere, can we draw a geographical map, with a given exponent, such that these three points are represented by three arbitrarily chosen points on the map?} Lagrange, in his paper, says that a geometric solution seems very difficult, and that he did not try to find a solution using algebra. Darboux solves the problem in a geometric manner. He says that it is the recent progress in geometry that made this solution possible.
In the second paper, \cite{Darboux-Chebyshev}, Darboux gives a detailed proof of the result of Chebyshev that we quoted above, saying that the most advantageous representation of a region of the sphere onto the Euclidean plane is the one where the magnification ratio is constant on the boundary of the surface to be represented. The proof he gives is the one outlined by Chebyshev, using potential theory. To make the relation with modern works, let us mention that Milnor also gave a proof of the theorem of Chebyshev, again following the latter's approach, see \cite{Milnor}. It seems that Milnor was not aware of the work of Darboux on the same problem (he does  not mention his name in his paper).

 We end the discussion on geography by recommending to the interested reader the book  \cite{Feeman} by T. G. Feeman.  It is addressed to mathematicians. It gives a nice exposition of the various distortions of a geographical map (angle, distance and area-distortions), with a review of several types of maps and an analysis of their distortion. The book also contains a section on the Tissot indicatrix.

 \section{Quasiconformal mappings}

 The theory of quasiconformal mappings appeared in complex analysis  under various names, and with slighly different definitions. 
  Presumably, the English name ``quasiconformal mapping" is due to Ahlfors; cf. his  comments in his \emph{Collected works edition}.\footnote{Ahlfors writes in \cite{Ahlfors-Collected}, Vol. 1, pp. 213: ``The truth is that I cannot recollect having invented the name, but I have also not been able to locate it elsewhere. Little did I know at the time what an important role quasiconformal mappings would come to play in my own work." In any case, the expression appears in print in Ahlfors' 1935  fundamental paper on covering surfaces \cite{Ahlfors1935}.}  Gr\"otzsch, in his papers written at the end of the 1920s, used the expression ``nichtkonformen Abbildungen" (non-conformal maps). He considered in \cite{Groetzsch1929} and in other papers he wrote in the same period the problem of finding the homeomorphism between two rectangles (the images of the vertices being fixed) that has the least deviation from conformality. He proved that in some reasonable sense of the expression `` that has the least deviation from conformality" the solution is given by the natural linear map between these rectangles. According to Ahlfors \cite{Ahlfors1964}, the problem that Gr\"otzsch solved in that paper was first considered as a mere curiosity, and the full strength of quasiconformal mappings and their use in the deformation theory of Riemann surfaces was first realized by Teichm\"uller. Let us quote Ahlfors from his 1954 paper  \cite{Ahlfors1964} pp. 156: 
\begin{quote} \small 
The very genesis of quasiconformal mappings was connected with the elementary extremal problem formulated by Gr\"otzsch. Teichm\"uller was the first to extract a general principle: In a class of mappings it is required to find one whose maximal dilatation is a minimum. It is to be expected that the solution is unique, and that the extremal mapping is characterized by simple properties.
\end{quote}
  
 Ahlfors gives a brief summary of the early use of quasiconformal mappings in his 1964 paper \cite{Ahlfors1964} (p. 153). In his 1978 ICM paper, he writes (\cite{Ahlfors1980} pp. 72): ``Quasiconformal mappings might have remained a rather obscure and peripheral object of study if it had not been for Oswald Teichm\"uller." 
  Quasiconformal mappings play a central role in the work of Teichm\"uller, who thoroughly  developed the theory and made it at the basis of several research topics, including the theory of moduli of Riemann surfaces, value distribution theory of meromorphic functions, the type problem and the Bieberbach coefficient problem. 
  We refer the reader to the surveys  \cite{2012d}, \cite{2015a} and \cite{Papa-qc} for more details.

  Lavrentieff, around 1935, wrote two papers in French, \cite{Lavrentiev1935a} and \cite{Lavrentiev1935}, in which he introduced a class of mappings he called ``fonctions presque analytiques" (almost analytic functions). The two papers  contain a remarkable series of results in function theory and geometry that are based on quasiconformal mappings. Let us note that in Lavrentieff's papers, the dilatation of a mapping between surfaces, in the modern sense of quasiconformal theory, is highlighted, but it is not assumed to be uniformly bounded. We note by the way that the same holds in some of Teichm\"uller's papers, e.g. \cite{T200}.

 We mention some of the results of Lavrentieff, because they remain poorly known compared to those of the other founders of the theory of quasiconformal mappings, and because the notion (which we recall below) of quasiconformal mappings that he uses is very close in spirit to the work of Tissot.

  In his paper \cite{Lavrentiev1935} (with an announcement of the results in the \emph{Comptes Rendus} note \cite{Lavrentiev1935a}), Lavrentieff gave a series of extremely interesting results on quasiconformal mappings. One of them (\S 4 of \cite{Lavrentiev1935a} and \S 3 of \cite{Lavrentiev1935}) is a generalization of Picard's theorem to the setting of quasiconformal mappings.  Another one  (\S 5 of \cite{Lavrentiev1935a} and \S 4 of \cite{Lavrentiev1935})  concerns a criterion, based on quasiconformal mappings, to find the type of a Riemann surface. There are several other results. We reproduce here Lavrentieff's definition of quasiconfomality to show how close it is to the ideas of Tissot. 
  
  Lavrentieff  says that  a function $w=f(z)$ of a complex variable $z$ in a domain $D$ of the complex plane is \emph{almost analytic} if it satisfies the following properties:
  
  \begin{enumerate}
  \item $f$ is single-valued and continuous on $D$.
  
  \item Except for a countable set of points $z$ in $D$, the function $f$ is an orientation-preserving local homeomorphism.
    
  \item There exist two real functions $p(z)\geq 1$ and $\theta(z)$ such that
  \begin{itemize}
  \item With the exception of points $z$ in a set $E$
   consisting of a finite number of analytic arcs, $p(z)$ is continuous, and $\theta(z)$ is continuous at all points $z$ satisfying $p(z)\not=1$.
  
  \item In every domain $\Delta$ which does not contain points of $E$ and whose frontier is a simple analytic curve, $p(z)$ is uniformly continuous, and if $\delta$ and its frontier do not contain points $z$ satisfying $p(z)=1$, $\theta$ is uniformly continuous on $\Delta$.
  
  \item Let $z_0$ be a point in $\Delta$ which is not in $E$, and consider the ellipse $\mathcal{E}$ centered at $z_0$, such that the angle between the major axis of $\mathcal{E}$  and the real axis of the plane is $\theta$, and such that if $a$ and $b$ are the major and the minor axes of $\mathcal{E}$, we have $p(z_0)=a/b$. 
  Then, we have
  \[\lim_{a\to 0}\big\vert\frac{f(z_1)-f(z_0)}{f(z_2)-f(z_0)} \big\vert =1
  \]
  where $z_1$ and $z_2$ are the points in $\mathcal{E}$ such that $\vert f(z)-f(z_0)\vert$ attains respectively its maximum and its minimum.
  \end{itemize}
  \end{enumerate}
  Lavrentieff calls the functions $p(z)$ and $\theta(z)$ the \emph{characteristic functions} of the almost analytic function $f(z)$. Although he uses the same word as Tissot, it seems that Lavrentieff was not aware of his work. 
  
  Lavrentieff notes that the problem of representing a 2-dimensional Riemannian manifold on a domain in the Euclidean plane is equivalent to the problem of constructing an almost-analytic function having the given characteristics $p$ and $\theta$. He gives several existence theorems that answer this question. We state his Theorem 3, pp. 414 of  \cite{Lavrentiev1935}:
  
  \medskip
  
  \emph{For any function $p(z)\geq 1$ and $\theta(z)$ defined on the unit disc $z\leq 1$ satisfying the conditions stated in the definition of an almost analytic function, one can construct an almost analytic function $w=f(z)$ satisfying $f(0)=0$, $f(1)=1$ which realizes a conformal representation from the unit disc $z\leq 1$ and the unit disc $w\leq 1$ and which has $p$ and $\theta$ as characteristic functions. 
  }
 
   \medskip
  
  This is one form of the fact which says (in modern terms) that an almost complex structure on a surface is integrable, which in turn is related to Gauss's existence of isothermal coordinates.
  
  In a note (p. 408), Lavrentieff says that in the definition of almost analytic functions if one adds the requirement that the characteristic function $p$ is bounded, then one obtains a class of functions which coincides with the ones considered by Gr\"otzsch in \cite{Gr1928}  (1928). He also notes that he had already considered a special class of almost analytic functions in his ICM paper \cite{Lavrentiev-ICM} (1928), in which he considers the question of constructing the Riemann Mapping Theorem by a sequence of explicit mappings obtained from the theory of partial differential equations, using a minimization principle. A similar application of quasiconformal mappings is mentioned by Teichm\"uller in the last part of his paper \cite{T20}. 
  We refer the reader to the recent survey \cite{Papa-qc}.

   \medskip
   
   \noindent {\bf Acknowledgements.---} I would like to thank Norbert A'Campo, Jeremy Gray, and Fran\c cois Laudenbach who read a preliminary version of this paper. The author acknowledges partial support of the French ANR (Agence Nationale de la Recherche) under the program ``G\'eom\'etrie de Finsler et applications."

 \end{document}